%%%%%%%%%%%%%%
%
%  Final draft
%
%  June 2, 2011
%
%%%%%%%%%%%%%%

\magnification=\magstephalf
\input amstex
\documentstyle{amsppt}
\pageheight{9truein}
\pagewidth{6.5truein}
\NoBlackBoxes
\TagsOnRight

\define \fa {\frak A}
\define \fb {\frak B}
\define \fc {\frak C}
\define \fd {\frak D}
\define \fI {\frak I}
\define \fp {\frak P}
\define \fw {\frak W}

\define \ord {\text{\rm ord}_v}
\define \ordw {\text{\rm ord}_w}
\define \divi {\text{\rm div}}
\define \La {\Lambda}
\define \fma {{\frak m}(A)}
\define \twedge {{\textstyle{\bigwedge}}}

\define \ba {\Bbb A}
\define \br {\Bbb R}
\define \bc {\Bbb C}
\define \be {\Bbb E}
\define \bq {\Bbb Q}
\define \qbar {\overline{\bq}}
\define \bp {\Bbb P}
\define \bfF {\Bbb F}
\define \bfq {\bfF _q}
\define \bfp {\bfF _p}
\define \fptbar {\overline{\bfp (T)}}
\define \fpbar {\overline{\bfp}}
\define \bz {\Bbb Z}

\define \ka {K_{\ba}}
\define \kbar {\overline{k}}
\define \nbar {\overline {N}}
\define \zetak {\zeta _K}

\define \ux {\bold x}
\define \us {\bold s}
\define \ue {\bold e}
\define \uv {\bold v}
\define \ua {\bold {\alpha}}
\define \ubeta {\bold {\beta}}
\define \udelta {\bold {\delta}}
\define \ugamma {\bold {\gamma}}
\define \ub {\bold b}
\define \uo {\bold 0}
\define \uy {\bold y}
\define \uY {\bold Y}
\define \uX {\bold X}
\define \uL {\bold L}
\define \uV {\bold V}
\define \ur {\bold r}
\define \uw {\bold w}

\define \nksep {N_k^{\text{sep}}}
\define \nfksep {NF_k^{\text{sep}}}

\topmatter
\title {Counting points of fixed degree and given height over function fields}
\endtitle
\author Jeffrey Lin Thunder and Martin Widmer
\endauthor
\thanks Second author supported by Swiss National Science Foundation
Grant \#118647 and Austrian Science Fund Grant \#M1222\endthanks
\address Dept. of Mathematics, Northern Illinois University, DeKalb, IL 60115,
USA
\endaddress
\email jthunder\@ math.niu.edu\endemail
\address Dept. for Analysis and Computational Number Theory, Graz Univ. of
Technology, Steyrergasse 30/II, A-8010 Graz, Austria
\endaddress
\email widmer\@ math.tugraz.at\endemail


\endtopmatter
\document
\baselineskip=20pt
\head Introduction\endhead

The purpose of a height in Diophantine geometry is to give a quantitative
measure of the arithmetic complexity of a point on some variety. This has
become a very important tool. Given a variety, one would like to know if
there are only finitely many points of a given height or height less than
a given bound. If so, one would further like to know the number of such points,
or at least upper and lower bounds for the number of such points. Another
goal would be to find asymptotic estimates for the number of such points as
the bound tends to infinity. Before we can discuss particular results here,
we must first set some notation.

For any field $k$ and a point $P=(x_1\:\cdots\: x_n)\in
\bp ^{n-1}(\overline{k})$ in
projective $n-1$-space over an algebraic closure,
let $k(P)$ denote
the field of definition of $P$ over $k$; in other words,
$k(P)$ is the field obtained by adjoining all possible quotients
$x_i/x_j$ to the field $k$. For a number field $k$, integers $n,\ d$ and
positive real $B$, let
$N_k(n,d,B)$ denote the number of points
in projective space
$\bp ^{n-1} (\qbar )$
with height less than $B$ and $[k(P)\: k]\le d$. An important early result
in this area is due to Northcott, who proved that
$N_k(n,d,B)$ is finite. On the other hand,
one easily sees that $N_k(n,d,B)$ grows without bound as $B\rightarrow
\infty$. Thus, one can ask for an asymptotic estimate, and it is precisely such
estimates that interest us here.

With this notation, $N_{\bq}(n,1,B)$ is simply the number of primitive
lattice points in a ball or cube (depending on the exact height used);
asymptotic results for $N_{\bq}(n,1,B)$ as $B\rightarrow\infty$
are classical. More generally, for any number
field $k$
Schanuel [Scha] proved that
$$N_k(n,1,B)=S_k(n,1)B^{ne}+O(B^{ne -1})$$
as $B\rightarrow\infty$,
where $e=[k\: \bq ]$ and the implicit constant depends on the number field
$k$ and the dimension $n$. (A $\log B$ term must be inserted in the error
when $k=\bq$ and $n=2$, though this is not necessary if one uses $L_2$ norms
at the infinite place.)
Here $S_k(n,1)$ is an explicitly given constant depending on the field $k$ and
the dimension $n$.
It turns out that proving similar asymptotic results for $N_k(n,d,B)$ when
$d>1$ is much more difficult.
Schmidt in [Schm1] gave non-trivial upper and lower bounds for
$N_{k}(n,d,B)$ and later gave asymptotic estimates for $N_{\bq }(n,2,B)$
in [Schm2].
In an unpublished thesis
Schmidt's student Gao found an asymptotic result
for $N_{\bq}(n,d,B)$ when $n>d+1>3$. Masser and Vaaler proved an asymptotic
result for $N_k(2,d,B)$ when $k=\bq$ in [MV1], and then for arbitrary
number fields in [MV2].  More recently,
the second author in [W1] proved that
$$N_k(n,d,B)=S_k(n,d)B^{ned}+O(B^{ned-1})$$
under the assumption that $n>5d/2+3+2/(ed)$, where $e=[k\:\bq ]$
as above and $S_k(n,d)$
is the sum of Schanuel constants $S_K(n,1)$ over extension fields $K$
of degree $d$ over $k$.
Our goal here is
to prove a result analogous to the theorem of Widmer above where  the number
field $k$
is replaced by a function field
over a finite field.

Fix a prime $p$, let $\bfp$ denote the finite field with $p$ elements and let
$T$ be transcendental
over this field, so that $\bfp (T)$ is a field of rational functions.
Throughout this article we fix algebraic closures $\fpbar$ of
$\bfp$ and $\fptbar\supset\fpbar$ of
$\bfp (T)$.
By a {\it function field} we will mean a finite
algebraic extension field $k\supseteq
\bfp (T)$ contained in $\fptbar$.
For such a field $k$ we have $k\cap\fpbar =\bfF _{q_k}$ for some finite field
$\bfF _{q_k}$; this is called the {\it field of constants} for $k$. Further,
we will write
$g_k$ for the genus, $J_k$ for the number of divisor classes of degree 0
(this is also the cardinality of the Jacobian) and $\zeta _k$ for the
usual zeta function of $k$. Serre stated, and later (independently) Wan [Wa] and
DiPippo [D] proved an analog for Schanuel's result in this context, where
the ``Schanuel constant" for a function field $k$ is
$$S_k(n,1)={J_k\over (q_k-1)\zeta _k(n)q_k^{n(g_k-1)}}.\tag 1$$

We denote the absolute additive height on $\fptbar$ by $h$ (defined
below). Fix a function field $k$ and set $e=[k\: \bfF _{q_k}(T)]$.
Suppose $K\supseteq k$ is another function field with $q_K=q_k$. Under
these assumptions,
the {\it effective degree} (see [A, chap. 15, \S 1])
of the extension field $K$ over $k$ is just the degree
$d=[K\: k]$, and the effective degree of $K$ over $\bfp (T)$ is $ed$. In
this case
the height of a point $P$ with $k(P)=K$ is
necessarily of the form $h(P)=m/ed$,
where $m$ is a non-negative integer.
This is a major difference from the number field situation, and leads
us to count not points of height no greater than a given bound, but
equal to a possible given bound. Moreover, since the possible heights
are naturally indexed by the non-zero integers, we are lead to the following
counting function.
\proclaim{Definition} Let $k$ be a function field and set
$e=[k\: \bfF _{q_k}(T)]$.
For integers $n>1$, $d\ge 1$ and $m\ge 0$,
$N_k(n,d,m)$ denotes the number
of points $P\in\bp ^{n-1}(\fptbar )$ with height $h(P)=m/ed$ and
$k(P)=K$ for some function field $K$ of degree $d$ over $k$ with $q_K=q_k$.
\endproclaim

Our main result is the following.

\proclaim {Theorem 1} Fix a function field $k$.
For all integers $n$ and $d>1$ satisfying
$n>d+2$, the sum
$$S_k(n,d)=\sum \Sb [K\: k]=d\\ q_K=q_k\endSb S_K(n,1)$$
converges.
Moreover, if $n>2d+3$ and $\varepsilon >0$ with $n>2d+3+\varepsilon$, then
for all
integers $m\ge 0$ we have
$$N_k(n,d,m)=S_k(n,d)q_k^{mn}+O\big (
q_k^{{m\over 2}(n+2d+3+\varepsilon )}\big ),$$
where the implicit constant depends only on $k$, $n$, $d$ and $\varepsilon$.
\endproclaim

The case where $d=1$ is just the function field version of Schanuel's
theorem. In that case one can do much better; see Theorem 2 below. We
will show that an asymptotic estimate for $N_k(n,d,m)$ of the form given
in Theorem 1 can only be possible when $n\ge d+1$.

Though Theorem 1 counts those points generating an extension of degree $d$ and
effective degree $d$, it is a simple matter to estimate the number of points
of given height generating an extension of degree $d$ and effective degree
$d'$ (which necessarily is a divisor of $d$) once we have Theorem 1.
We note that
the height of such a point is necessarily of
the form $m/ed'$ for some non-negative integer $m$.
Let $k,\ d$ and $e$ be as in the statement of Theorem 1
and suppose one wants to count the number $N$ of points $P$ with
$[k(P)\: k]=d$, $q_{k(P)}=q_k^{d/d'}$ and $h(P)=m/ed'$.
Certainly all such points will be counted in
$N_{k\bfF _{q_k^{d/d'}}}(n,d',m)$, where $k\bfF _{q_k^{d/d'}}$ denotes
the compositum field, but this will be an over-count since
we have $[k(P)\:k]\le [k\bfF _{q_k^{d/d'}}(P)\: k\bfF _{q_k^{d/d'}}]\cdot
[\bfF _{q_k^{d/d'}}\: \bfF _{q_k}]$. Put another way, if $P$ is counted
in $N_{k\bfF _{q_k^{d/d'}}}(n,d',m)$, then we have $q_{k(P)}=q_k^r$ for
some $r\le d/d'$.
We thus see that
$$N_{k\bfF _{q_k^{d/d'}}}(n,d',m)\ge N\ge
N_{k\bfF _{q_k^{d/d'}}}(n,d',m)-\sum _{1\le r <d/d'}
N_{k\bfF _{q_k^r}}(n,d',m).$$
Note that the summands subtracted here are of a lower order of magnitude
than $N_{k\bfF _{q_k^{d/d'}}}(n,d',m)$ by Theorem 1, whence we have an
asymptotic estimate for the desired quantity $N$.

We can also use Theorem 1 to count certain forms.
Suppose $F(\uX )\in k[\uX]$ is a homogeneous polynomial (form) in
$n$ variables of degree $d$. Such a form is called {\it decomposable} if
it factors completely into a product of $d$ linear forms:
$$F(\uX )=\prod _{i=1}^dL_i(\uX ).$$
Denote the coefficient vector of the linear factor $L_i(\uX)$ by
$\uL _i$. Clearly these $\uL _i$ are unique only up to a scalar multiple;
we thus identify each $\uL _i$ with a point $P(\uL _i)\in\bp ^{n-1}(\fptbar )$.
In a similar manner, we identify the set of proportional forms
$\lambda F(\uX )$ with a point $P(F)\in\bp ^{{d+n-1\choose n-1}-1}(k)$. Thus
the number of non-proportional forms $F(\uX )\in k[\uX ]$ of
degree $d$ in $n$ variables with $h\big (P(F)\big )=m/e$ is exactly
$N_k\big ({d+n-1\choose n-1},1,m\big )$. We can also use Theorem 1 to
count certain decomposable forms.

\proclaim{Definition} Let $k$ be a function field and set
$e=[k\: \bfF _{q_k}(T)]$.
Fix positive integers $n$ and $d$, and a non-negative integer $m$.
Then $NF_k(n,d,m)$ denotes the number of non-proportional decomposable forms
$F(\uX )\in k[\uX ]$  in $n$ variables of degree $d$
with height $h\big (P(F)\big )=m/e$, where each $k\big (P(\uL _i)\big )$
is an extension of degree $d$ and effective degree $d$ over $k$.
\endproclaim

As we noted above, the height $h\big (P(F)\big )$ is necessarily of the form
$m/e$ for some integer $m$ whenever $F(\uX )\in k[\uX ]$.

\proclaim{Corollary} Fix a function field $k$ and positive integers
$n$ and $d$. Let $p^r$ denote the highest power of $p$ dividing $d$,
where $p$ is the characteristic of $k$. Then for all
integers $m\ge 0$ we have
$$dNF_k(n,d,m)=N_k(n,d,m)+\sum _{i=1}^r(p^i-p^{i-1})N_k(n,d/p^i,m).$$
(As usual, empty sums are to be interpreted as zero.)
In particular, if $n>2d+3$ and $\varepsilon >0$ with $n>2d+3+\varepsilon$, then
for all integers $m\ge 0$ we have
$$\aligned dNF_k(n,d,m)&=q_k^{nm}
\left (S_k(n,d)+\sum _{i=1}^r (p^i-p^{i-1})S_k(n,d/p^i)\right )\\
&\qquad +
O\big (
q_k^{{m\over 2}(n+2d+3+\varepsilon)}\big ),\endaligned$$
where the implicit constant depends only on $k$, $n$, $d$ and $\varepsilon$.
\endproclaim

We will give a proof of the Corollary after our proof of Theorem 1 in
the final section.
We conclude our introduction with a bit more notation and the definition
of the height used above. In the next section we outline our method of proof
and state its main ingredients. The following sections are devoted to
auxiliary results and the proofs of our main theorems and their corollaries.

For a function field $k$ let
$M(k)$ denote the set of places of $k$. For every place $v\in M(k)$
let $k_v$ denote the topological
completion of $k$ and let $\ord$ denote the order function on $k_v$ normalized
to have image $\bz \cup \{\infty\}.$
We extend $\ord$ to $k_v^n$ by defining
$$\ord(x_1,...,x_n)=\min_{1\leq i\leq n}\{\ord x_i\},$$
with the usual convention that $\ord 0=\infty$ (greater than any integer).
Each non-zero element $\ux$ of $k^n$ gives rise to a divisor
$$\divi (\ux )=\sum_{v\in M(k)}\ord (\ux)\cdot v.$$
For such an $\ux$ we define the {\it relative height} to be
$$h_k(\ux )=-\deg\divi (\ux ).$$
Clearly $h_k$ is an integer. Moreover, since the degree of a principal
divisor is 0, $h_k$ is actually a function on projective space.
In particular, we can assume without loss of generality that one
of the coordinates of $\ux$ is 1, so that $\ord (\ux )\le 0$ for
all places $v$ and $h_k(\ux )$ is necessarily a
non-negative integer.
Now $[k\: \bfF _{q_k} (T)]$ is, by definition,
the effective degree of
the extension $k$ over $\bfp (T)$; we define the {\it absolute height} $h$
to be
$$h(\ux )={h_k(\ux )\over [k\: \bfF _{q_k}(T)]}.$$

Dividing the relative height by the effective degree gives a height
that is not dependent on the choice of field. Specifically, if
$P\in\bp ^{n-1}(\fptbar )$ is defined over $k$ and $K$ is any function field
containing $k$, so that $P$ is in both $\bp ^{n-1}(k)$ and $\bp ^{n-1}(K)$,
we have
$$h_K(P)=h_k(P){[K\: k]\over [\bfF _{q_K}\: \bfF _{q_k}]}
=h_k(P){[K\: \bfF _{q_K}(T)]\over [k\: \bfF _{q_k}(T)]}.$$
(See [T1, p. 150]). Thus, the height $h$ is a
function on $\bp ^{n-1}(\fptbar )$.

\head Outline of the proof\endhead

As one would suppose from the statement of Theorem 1, we estimate
$N_k(n,d,m)$ by summing over all possible function fields $K\supseteq k$
of degree $d$ with $q_K=q_k$. More precisely, for such a field $K$
we let $N_k(n,K,m)$ denote the number of points $P\in\bp ^{n-1}(K)$ with
$h(P)=m/de$ and $k(P)=K$, where $e=[k\:\bfF _{q_k}(T)]$. In other words,
$N_k(n,K,m)$ is the number of those points $P$ counted in $N_k(n,d,m)$
where $k(P)=K$ for the fixed field $K$. We thus have
$$N_k(n,d,m)=\sum \Sb [K\: k]=d\\ q_K=q_k\endSb N_k(n,K,m).\tag 2$$

It will not be difficult to prove that the main term in the
asymptotic estimate for
$N_k(n,K,m)$ is $S_K(n,1)q_K^{nm}$. Our efforts will mainly be focused
on the error term.
The major ingredient of our
proof is a version of Schanuel's result
for function fields where we pay particular attention to the form of the
error term.

\proclaim{Theorem 2} Let $k$ be a function field
and set $e=[k\:\bfF _{q_k} (T)]$. Suppose $m$ is an integer with
$m\ge 2g_k-1$ and $1/4\ge\varepsilon >0$. Then
for all integers $n\ge 4$ we have
$$N_k(n,1,m)=S_k(n,1)q_k^{nm}+O\big (q_k^{m(1+\varepsilon )}
q_k^{g_k(n-2-2\varepsilon)}\big ),$$
and for $n=2,3$
$$N_k(n,1,m)=S_k(n,1)q_k^{nm}+O\big (q_k^{m(1+\varepsilon )}
q_k^{g_k(1+\varepsilon )}\big ).$$
Suppose $m< 2g_k-1$. Then
for all $\varepsilon >0$ and all integers $n\ge 2$ we have
$$N_k(n,1,m)\ll
q_k^{m({n+1\over 2}+\varepsilon )}$$
All the implicit constants here depend only on $n$, $e$,
$q_k$ and $\varepsilon$.
\endproclaim

\proclaim{Corollary 1} Let $k$ be a function field. For all integers $n\ge 2$
and $m\ge 0$
$$N_k(n,1,m)\ll q_k^{nm},$$
where the implicit constant depends only on $n$ and $q_k$.
\endproclaim

We will also use the following quantity.

\proclaim{Definition} For function fields $K\supseteq k$ and integers
$n>1$,
$$\delta _n(K/k)=\min\{ h(P)\: P\in\bp ^{n-1}(K),\ k(P)=K\}.$$
\endproclaim

\proclaim{Corollary 2}
Let $K\supseteq k$ be function fields with $q_K=q_k$, set $[K\: k]=d$ and
$e=[k\: \bfF _{q_k}(T)]$. Suppose $m$ and $n$ are positive integers satisfying
$m\ge de\delta _n(K/k)$ and $\varepsilon >0$. Then if $n\ge 4$
we have
$$N_k(n,K,m)=S_K(n,1)q_k^{nm}+\cases
O\big (q_k^{nm/2}\big )&\text{if $m\ge 2g_K-1$,}\\
O\big (q_k^{m({n+1\over 2}+\varepsilon )}\big )&
\text{otherwise,}\endcases$$
and if $n=2,3$
$$N_k(n,K,m)=S_K(n,1)q_k^{nm}+\cases
O\big (q_k^{{3m\over 2}(1+\varepsilon )}\big )&\text{if $m\ge 2g_K-1$,}\\
O\big (q_k^{m({n+1\over 2}+\varepsilon )}\big )&
\text{otherwise.}\endcases$$
The implicit constants here depend only on $n$, $e$, $d$, $q_k$ and
$\varepsilon$. For all integers $m<de\delta _n(K/k)$, $N_k(n,K,m)=0$ by
definition.
\endproclaim

We will see that the ``main terms" here are majorized by the ``error terms"
in the case $m< 2g_K-1$. It will be convenient for our purposes to have
a uniform statement, however.
Our proof of Theorem 1 will
use Corollary 2 and (2).

\head Proof of Theorem 2 and its Corollaries\endhead

Our proof of Theorem 2 will follow along the same
lines as  the proof of Theorem 1 of [T2].
Our job here is made easier since we don't look at arbitrary ``twisted"
heights, but we need to work somewhat harder to get good explicit
dependencies on the field. Throughout this section
all function fields appearing are assumed
to have the same field of constants; we will
write $q$ for the cardinality of this field. Before we get to the proof
of Theorem 2, we need to recall some
concepts from the theory of function fields and
prove a few auxiliary results. In what follows, divisors will always be
denoted using capital script German font ($\fa$, $\fb$, etc.),
with the exception of the zero divisor which will be
denoted by $0$.

Let $k$ be a function field and $n$ be a positive integer. For a
divisor $\fa$, set
$$L(\fa ,n)=\{ \ux\in k^n\: \ord (\ux )\ge -\ord (\fa )\ \text{for all
$v\in M(k)$}\}.$$
Then $L(\fa ,n)$ is a vector space of finite dimension over $\bfq$
(see [T1,\S II]); we denote this dimension by $l(\fa ,n)$. Thus,
the cardinality of $L(\fa ,n)$ is $q^{l(\fa ,n)}.$ It will prove convenient
to write $\lambda (\fa ,n)$ for the number of non-zero elements of $L(\fa ,n)$,
i.e., $\lambda (\fa ,n)=q^{l(\fa ,n)}-1.$ Let $L'(\fa ,n)$ denote the
set of
those $\ux\in L(\fa ,n)$ with $\ord (\ux )=-\ord (\fa )$ for all places
$v\in M(k)$ and write $\lambda '(\fa ,n)$ for its cardinality.

\proclaim{Lemma 1} For a function field $k$ and divisor $\fa$ we have
$$l(\fa ,n)=nl(\fa ,1)=n\left (\deg (\fa )+1-g_k+l(\fw -\fa ,1)\right ),$$
where $\fw$ is any divisor in the canonical class. In particular,
$l(\fa ,n)=n\big (\deg (\fa )+1-g_k\big )$
whenever $\deg (\fa )\ge 2g_k-1$,
$l(\fa ,n)\le {n\over 2}\big (\deg (\fa )+1\big )$
whenever $0\le\deg (\fa )\le 2g_k-2$, and
$l(\fa ,n)=0=\lambda (\fa ,n)$ whenever $\deg (\fa )<0$.
\endproclaim

\demo{Proof}
One readily sees that $l(\fa ,n)=nl(\fa ,1)$. The lemma thus follows from
the Riemann-Roch Theorem and Clifford's Theorem (see [S, Chap. 1], for
example).
\enddemo

Next, for all integers $l\ge 0$ write $a(l)$ for the number of non-negative
divisors of degree $l$:
$$a(l)=\sum \Sb \fc\ge 0\\ \deg (\fc )=l\endSb 1.$$
Then the zeta function is given by
$$\sum _{l=0}^{\infty}a(l)q^{-sl}=\zeta _k(s)$$
for all $s>1$. We let $\mu$ denote the usual
M\"obius function on the divisor group. It is defined by the
following four conditions: $\mu (0)=1$, $\mu (\fa+\fb )=\mu (\fa )\mu (\fb )$
whenever $\fa$ and $\fb$ are relatively prime (i.e., have disjoint support),
$\mu (\fp )=-1$ if $\fp$ is a prime divisor, and $\mu (r\fp )=0$ if
$\fp$ is a prime divisor and $r>1$. Write
$$b(l)=\sum \Sb \fc \ge 0\\ \deg (\fc )=l\endSb \mu (\fc ).$$
Then as is well-known (see [T2, Lemma 4], for example)
$$\sum _{l=0}^{\infty}b(l)q^{-sl}={1\over\zeta _k(s)}\tag 3$$
for all $s>1$.

\proclaim{Lemma 2}
Fix a function field $k$ and set $e=[k\: \bfq (T)]$.
Then  $1<\zeta _k(s)\le \big (\zeta _{\bfq (T)}(s)\big )^e$
for all $s>1$.
For all integers $m\ge 0$, all $s\le 1$ and all $\varepsilon >0$ we have
$$\sum _{l=0}^ma(l)q^{-sl}\ll q^{m(1-s+\varepsilon)},$$
and  for all $s>1+\varepsilon$
$$\sum _{l\ge m}a(l)q^{-sl}\ll q^{-m(s-1-\varepsilon )},$$
where the implicit constants depend only on $q, \ e$ and $\varepsilon$.
In particular,
$$a(m)\ll q^{m(1+\varepsilon )}$$
for all integers $m\ge 0$ and all $\varepsilon >0$. Finally,
$a(m)={J_k\over q-1}(q^{m+1-g_k}-1)$
for all integers $m\ge 2g_k-1.$
\endproclaim

\demo{Proof} We have (see [S, V.1.4 Lemma], for example)
$$a(m)={1\over q-1}\sum _{j=1}^{J_k}q^{l(\fc _j ,1)}-1,\tag 4$$
where $\fc _1,\ldots ,\fc _{J_K}$ are representatives of the divisor classes
of degree $m$. In particular, $a(0)=1$ and
$a(m)={J_k\over q-1}(q^{m+1-g_k}-1)$ for all $m\ge 2g_k-1$ by the Riemann-Roch
Theorem. We get $1<\zeta _k(s)$ at once.
Since the genus is 0 and the number of divisor classes of
degree 0 is 1 for a field of rational functions, we get the well-known formula
$$\zeta _{\bfq (T)}(s)={1\over (1-q^{-s})(1-q^{1-s})}$$
for all $s>1$. In particular,
$$\zeta _{\bfq (T)}(1+\varepsilon )\ll 1.$$

We next use the Euler product
(see [S, V.1.8 Proposition], for example):
$$\zeta _k(s)=\prod _{v\in M(k)}(1-q^{-s\deg (v)})^{-1}.$$
For a place $v\in M(k)$ lying over a place $w\in M(\bfq (T))$, write
$f_v$ for the residue class degree and $e_v$ for the ramification index.
Then as is well-known, $\sum _{v|w}e_vf_v=[k\:\bfq (T)]$ for all places
$w\in M(\bfq (T))$. Since the ramification indices $e_v$ are always positive
integers, we get
$$\aligned \prod _{v|w}(1-q^{-s\deg (v)})^{-1}&=\prod _{v|w}
(1-q^{-sf_v\deg (w)})^{-1}\\
&\le \big ((1-q^{-s\deg (w)})^{-1}\big )^{\sum _{v|w}f_v}\\
&\le \big ((1-q^{-s\deg (w)})^{-1}\big )^{\sum _{v|w}e_vf_v}\\
&= \big ((1-q^{-s\deg (w)})^{-1}\big )^{[k\:\bfq (T)]}\endaligned$$
for all places $w\in M(\bfq (T))$. Thus $\zeta _k(s)\le\big (
\zeta _{\bfq (T)}(s)\big )^e$.

Now if $s\le 1$ we have
$$\aligned
\sum _{l=0}^ma(l)q^{-sl}&\le\sum _{l=0}^mq^{(m-l)(1-s+\varepsilon )}a(l)
q^{-sl}\\
&=q^{m(1-s+\varepsilon )}\sum _{l=0}^ma(l)q^{-l(1+\varepsilon )}\\
&<q^{m(1-s+\varepsilon )}\sum _{l=0}^{\infty }a(l)q^{-l(1+\varepsilon )}\\
&=q^{m(1-s+\varepsilon )}\zeta _k(1+\varepsilon )\\
&\le q^{m(1-s+\varepsilon )}\big (\zeta _{\bfq (T)}(1+\varepsilon )\big ) ^e\\
&\ll q^{m(1-s+\varepsilon )},\endaligned$$
and if $s>1+\varepsilon$
$$\aligned
\sum _{l\ge m}a(l)q^{-sl}&=\sum _{l\ge m}a(l)q^{-(1+\varepsilon )l}q^{-(s-1-
\varepsilon )l}\\
&<q^{-m(s-1-\varepsilon )}\sum _{l\ge m}a(l)q^{-(1+\varepsilon )l}\\
&\le q^{-m(s-1-\varepsilon )}\zeta _k(1+\varepsilon )\\
&\ll q^{-m(s-1-\varepsilon )}.\endaligned$$

\enddemo

\proclaim{Lemma 3}
Fix a function field $k$.
Then for all $\varepsilon >0$
$$q^{g_k(1-\varepsilon )}\ll J_k\ll q^{g_k(1+\varepsilon )},$$
where the implicit constants depend only on $q$, $[k\: \bfq (T)]$
and $\varepsilon$.
\endproclaim

\demo{Proof}
If the genus is 0, then $J_K=1$ and the statment is true,
so assume that $g_k\ge 1$.

By Lemma 2
$$\aligned
J_kq^{g_k}&\ll
{J_k\over q-1}(q^{g_k}-1)\\
&=a(2g_k-1)\\
&\ll q^{(2g_k-1)(1+\varepsilon /2)}\\
&< q^{g_k(2+ \varepsilon )},\endaligned$$
so that
$J_k\ll q^{g_k(1+\varepsilon )}.$

By Lemma 2,
$$\aligned
\zeta _k(s)&=\sum _{l=0}^{\infty}a(l)q^{-sl}\\
&={J_k\over q-1}\sum _{l=2g_k-1}^{\infty}(q^{l+1-g_k}-1)q^{-sl}+
\sum _{l=0}^{2g_k-2}a(l)q^{-sl}\\
&={J_kq^{s(1-2g_k)}\over q-1}\left ({q^{g_k}\over 1-q^{1-s}}-{1\over 1-q^{-s}}
\right )+
\sum _{l=0}^{2g_k-2}a(l)q^{-sl}.\endaligned$$
This identity is used to analytically continue the zeta function to the
complex plane, with simple poles at $s=0,1$. Further, by the ``Riemann
Hypothesis", i.e., Hasse-Weil Theorem, this analytic continuation has
exactly $2g_k$ zeros (counting multiplicity), all of which have real
part equal to 1/2. In particular, the analytically continued zeta function
is negative for all $1/2<s<1$. Hence, setting $s=1-\varepsilon$, we have
for all positive $\varepsilon <1/2$
$$
{J_kq^{(1-\varepsilon )(1-2g_k)}\over q-1}\left ({q^{g_k}\over 1-q^{
\varepsilon }}-{1\over 1-q^{\varepsilon -1}}
\right ) +\sum _{l=0}^{2g_k-2}a(l)q^{l(\varepsilon -1)} <0.$$
Since we are assuming $g_k\ge 1$, we have $\sum _{l=0}^{2g_k-2}a(l)q^{l(
\varepsilon -1)}\ge a(0)=1$. Thus
$${J_kq^{(1-\varepsilon )(1-2g_k)}\over q-1}\left ({q^{g_k}\over 1-q^{
\varepsilon }}-{1\over 1-q^{\varepsilon -1}}\right ) <-1.$$
Multiplying both sides by $(1-q^{\varepsilon })
(1-q^{\varepsilon -1})$ (which is negative), we get
$$\aligned (q^{\varepsilon }-1)(1-q^{\varepsilon -1})&<
{J_kq^{(1-\varepsilon )(1-2g_k)}\over q-1}\big ( q^{g_k}(1-q^{\varepsilon
-1})-
(1-q^{\varepsilon})\big )\\
&\le {J_kq^{(1-\varepsilon )(1-2g_k)}\over q-1}
(q^{g_k}-1)\\
&<{J_kq^{g_k(2\varepsilon -1)}q^{(1-\varepsilon )}\over q-1}.\endaligned$$
This shows the other inequality.
\enddemo

We note that there always exists a divisor of degree 1 [S, V.1.11 Corollary].

\proclaim{Lemma 4} Let $k$ be a function field and suppose $n\ge 2$ is an
integer. Set representatives $\fa _1,\ldots ,\fa _{J_k}$ of the divisor
classes of degree 0 and fix a divisor $\fa _0$ of degree 1. Then
for all integers $0\le i\le 2g_k-2$ we have
$$\sum _{j=1}^{J_k}\lambda (\fa _j+i\fa _0,n)\ll
a(i)q^{(n-1)i/2},$$
where the implicit constant depends only on $n$ and $q$.
\endproclaim

\demo{Proof}
By (4) we have $\sum _{j=1}^{J_k}q^{l(\fa _j+i\fa _0,1)}-1=(q-1)a(i)$.
Setting $c_j=l(\fa _j+i\fa _0,1)$, we get $c_j\le {i+1\over 2}$ by
Lemma 1, whence
$$\aligned
\sum _{j=1}^{J_k}\lambda (\fa _j+i\fa _0,n)&=
\sum _{j=1}^{J_k}q^{nc_j}-1\\
&=
\sum _{j=1}^{J_k}\big (q^{c_j}-1\big )\big (q^{(n-1)c_j}+q^{(n-2)c_j}+\cdots
+1\big )\\
&\ll
\sum _{j=1}^{J_k}\big (q^{c_j}-1\big )q^{(n-1)i/2}\\
&\ll a(i)q^{(n-1)i/2}.\endaligned$$
\enddemo

\proclaim{Lemma 5}
Let $k$ be a function field and suppose $n\ge 2$ is an
integer. Set representatives $\fa _1,\ldots ,\fa _{J_k}$ of the divisor
classes of degree 0 and fix a divisor $\fa _0$ of degree 1.
For all integers $0\le i\le 2g_k-2$ we have
$$\sum _{j=1}^{J_k}\lambda (\fa _j+i\fa _0,n)-(q^{n(i+1-g_k)}-1)=
q^{n(i+1-g_k)}\sum _{j=1}^{J_k}\lambda (\fa _j+(2g_k-2-i)\fa _0,n).$$
\endproclaim

\demo{Proof}
Let $\fw$ be a divisor in the canonical class.
By the definition of $\lambda$ and Lemma 1,
$$\aligned
\lambda (\fa _j+i\fa _0,n)-(q^{n(i+1-g_k)}-1)&=q^{nl(\fa _j+i\fa _0,1)}-q^{n(
i+1-g_k)}\\
&=q^{n(i+1-g_k)}\left (q^{nl(\fw-\fa _j-i\fa _0,1)}-1\right )\endaligned$$
for all $0\le i\le 2g_k-2$ and $j=1,\ldots ,J_k$. Clearly
$\fw -\fa _j-i\fa _0$ runs through all divisor classes of degree $2g_k-2-i$
as $j$ goes from $1$ to $J_k$, since $\deg (\fw )=2g_k-2$. Thus
$$\sum _{j=1}^{J_k}\lambda (\fa _j+i\fa _0,n)-\big (q^{n(i+1-g_k)}-1\big )
=q^{n(i+1-g_k)}
\sum _{j=1}^{J_k}\lambda \big (\fa _j+(2g_k-2-i)\fa _0,n\big ).$$

\enddemo

\demo{Proof of Theorem 2}
To ease notation, write $J$ and $g$ for $J_k$ and $g_k$, respectively.
Set
representatives
$\fa _1,\ldots ,\fa _J$ of the divisor classes of degree 0 and
fix a divisor $\fa _0$ of degree 1. All implicit constants appearing
in our proof depend only on (at most) $n,$ $e$, $q$ and $\varepsilon$.

Using M\"obius inversion exactly as in [T2, \S 4], we get
$$\aligned
(q-1)N_k(n,1,m)&=\sum _{j=1}^J\lambda '(\fa _j+m\fa _0,n)\\
&=\sum _{j=1}^J\sum _{\fc\ge 0}\mu (\fc )\lambda (\fa _j+m\fa _0-\fc ,n)\\
&=\sum _{l=0}^m b(l)\sum _{j=1}^J\lambda (\fa _j+(m-l)\fa _0,n),
\endaligned\tag 5$$
where the last equation follows from the fact that $l(\fa ,n)=l(\fb ,n)$
whenever $\fa$ and $\fb$ are linearly equivalent divisors.

Now assume $m\ge 2g-1$. From (5) and Lemma 1 we have
$$\aligned (q-1)N_k(n,1,m)
&=\sum _{l=0}^m b(l)\sum _{j=1}^J\lambda (\fa _j+(m-l)\fa _0,n)\\
&=\sum _{j=1}^J\sum _{l=0}^{\infty}b(l)q^{n(m-l+1-g)}-
\sum _{j=1}^J\sum _{l=0}^{m}b(l) -
\sum _{j=1}^J\sum _{l=m+1}^{\infty}b(l)q^{n(m-l+1-g)}\\
&\qquad +
\sum _{j=1}^J\sum _{l=m-2g+2}^m
b(l)\left (\lambda (\fa _j+(m-l)\fa _0,n)-(q^{n(m-l+1-g)}-1)\right ).
\endaligned\tag 6$$

By (3)
$$\sum _{j=1}^J\sum _{l=0}^{\infty}b(l)q^{n(m-l+1-g)}={Jq^{n(m+1-g)}\over
\zeta _k(n)}.\tag 7$$
Clearly $a(l)\ge |b(l)|$ always, so that by Lemmas 2 and 3
$$
\left |\sum _{j=1}^J\sum _{l=0}^{m} b(l)\right |\le
\sum _{j=1}^J\sum _{l=0}^{m} a(l)
\ll Jq^{m(1+\varepsilon )}\ll q^{m(1+\varepsilon )}q^{g(1+\varepsilon )}\tag 8
$$
for all $\varepsilon >0$.
Similarly (and since $n\ge 2$)
$$\left |
\sum _{j=1}^J\sum _{l=m+1}^{\infty}b(l)q^{n(m-l+1-g)}\right |\le
\sum _{j=1}^J\sum _{l=m+1}^{\infty}a(l)q^{n(m-l+1-g)}
\ll Jq^{-ng}q^{m(1+\varepsilon )}\ll q^{m(1+\varepsilon )}.\tag 9$$

We now turn to the last term in (8). First, by Lemma 1 we have
$$0\le \lambda (\fa _j+(m-l)\fa _0,n)-(q^{n((m-l)+1-g)}-1)$$
for all $j=1,\ldots ,J$ and $l=m-2g+2,\ldots ,m$. Hence
$$\aligned &
\left |\sum _{j=1}^J\sum _{l=m-2g+2}^m
b(l)\left (\lambda (\fa _j+(m-l)\fa _0,n)-(q^{n(m-l+1-g)}-1)\right )\right |\\
&\qquad
\le
\sum _{j=1}^J\sum _{l=m-2g+2}^m
a(l)\left (\lambda (\fa _j+(m-l)\fa _0,n)-(q^{n(m-l+1-g)}-1)\right )\\
&\qquad\qquad
=\sum _{i=0}^{2g-2}\sum _{j=1}^Ja(m-i)
\left (\lambda (\fa _j+i\fa _0,n)-(q^{n(i+1-g)}-1)\right ),\endaligned\tag 10$$
where we have written $i$ for $m-l$.
Setting $i'=2g-2-i$, by Lemmas 2 (applied to $a(m+i'-2g+2)$), 4 and 5
$$\aligned
\sum _{i=0}^{2g-2}\sum _{j=1}^J&a(m-i)
\left (\lambda (\fa _j+i\fa _0,n)-(q^{n(i+1-g)}-1)\right )\\
&=\sum _{i'=0}^{2g-2}\sum _{j=1}^Ja(m+i'-2g+2)q^{n(g-1-i')}\lambda (
\fa _j+i'\fa _0,n)\\
&\ll\sum _{i'=0}^{2g-2}a(m+i'-2g+2)q^{n(g-1-i')}
a(i')q^{(n-1)i'/2}\\
&\ll q^{m(1+\varepsilon )}q^{g(n-2-2\varepsilon )}\sum _{i'=0}^{2g-2}a(i')
q^{i'(1+\varepsilon -(n+1)/2)}
.\endaligned \tag 11$$
If $n\ge 4$ and $\varepsilon \le 1/4$, then $(1+\varepsilon) -(n+1)/2\le -5/4$,
so that by (10), (11) and Lemma 2 
$$\multline
\left |\sum _{j=1}^J\sum _{l=m-2g+2}^m
b(l)\left (\lambda (\fa _j+(m-l)\fa _0,n)-(q^{n(m-l+1-g)}-1)\right )\right |\\
\ll q^{m(1+\varepsilon )}q^{g(n-2-2\varepsilon )}\sum _{i'=0}^{2g-2}a(i')
q^{i'(-5/4)}\\
< q^{m(1+\varepsilon )}q^{g(n-2-2\varepsilon )}\zeta _k(5/4)\\
\ll q^{m(1+\varepsilon )}q^{g(n-2-2\varepsilon )}.\endmultline\tag 12$$
If $n=2,3$  we use $a(i')\ll q^{i'(1+\varepsilon /2)}$, so by (10) and (11)
$$\aligned &
\left |\sum _{j=1}^J\sum _{l=m-2g+2}^m
b(l)\left (\lambda (\fa _j+(m-l)\fa _0,n)-(q^{n(m-l+1-g)}-1)\right )\right |\\
&\qquad\qquad\qquad\ll q^{m(1+\varepsilon )}q^{g(n-2-2\varepsilon )}
\sum _{i'=0}^{2g-2}a(i')q^{i'(1+\varepsilon -(n+1)/2)}\\
&\qquad\qquad\qquad\ll q^{m(1+\varepsilon )}q^{g(n-2-2\varepsilon )}
\sum _{i'=0}^{2g-2}q^{(i'/2)(4+3\varepsilon -(n+1))}\\
&\qquad\qquad\qquad\ll q^{m(1+\varepsilon )}q^{g(n-2-2\varepsilon )}q^{g(4+
3\varepsilon -(n+1))}\\
&\qquad\qquad\qquad =q^{m(1+\varepsilon )}q^{g(1+\varepsilon )}
.\endaligned\tag 13$$
The case where $m\ge 2g-1$ follows from (1),  (6)-(9), (12), (13)  and Lemma 2.

We now turn to the case where $m\le 2g-2$. By (5) and the definitions we have
$$\aligned (q-1)N_k(n,1,m)&=\sum _{j=1}^J\lambda ' (\fa _j+m\fa _0,n)\\
&\le \sum _{j=1}^J\lambda  (\fa _j+m\fa _0,n).\endaligned$$
The proof is completed by this and Lemmas 2 and 4.
\enddemo

\demo{Proof of Corollary 1}
By (1) and Lemmas 2 and 3,
$$S_k(n,1)<J_kq^{-ng_k}\ll 1\tag 14$$
for all function fields $k$ and
all integers $n\ge 2$. One readily verifies that
$$q^{m(1+\varepsilon )}q^{g_k(n-2-2\varepsilon )}\ll q^{nm/2}\tag 15$$
for all integers $n\ge 4,\
m\ge 2g_k-1$ and all $\varepsilon \le 1/4$. Also,
$$q^{m(1+\varepsilon )}q^{g_k(1+\varepsilon )}\ll q^{(3m/2)(1+\varepsilon )}
\tag 16$$
for all integers $m\ge 2g_k-1$ and all $\varepsilon >0$.
Corollary 1 follows from Theorem 2 and (14)-(16).
\enddemo

The proof of Corollary 2 will require one further auxiliary result.

\proclaim{Lemma 6}
Suppose
$K\supseteq k$ are function fields and write
$d=[K\: k]$. Then then number
$N$ of intermediate fields $L$ with $k\subseteq L\subseteq K$ satifies
$N\le d2^{d!}$.
\endproclaim

\demo{Proof}
Let $k\subseteq L\subseteq K$. Suppose first that $L$ is a separable
extension of $k$. We have $[L\: k]\le [K\: k]=d,$ whence by elementary Galois
theory at most $2^{d!}$ possible $L$. Now suppose that $L$ is
not a separable extension of $k$. Then we have $k\subseteq L_s\subset L$,
where $L_s$ is a separable extension of $k$ and $L$ is a purely
inseparable extension of $L_s$. Then we must have $[L\: L_s]=p^r$ for
some positive integer $r$ (recall that $p$ is the characteristic of all
our fields). Moreover, $L_s=\{ a^{p^r}\: a\in L\}$ (see [S, Proposition
III.9.2], for example).
Therefore each element in $L_s$ has a unique $p^r$-th root, so that $L_s$
and $p^r$ completely determine $L$. Since both $[L\: L_s]$ and $[L_s\: k]$
are no greater than $d$, we get our estimate.
\enddemo

\demo{Proof of Corollary 2}
Set $d=[K\: k]$. We then have
$$N_k(n,K,m)=N_K(n,1,m)-\sum \Sb d'<d\\ d'|d\endSb
\sum \Sb k\subseteq L\subset K\\
[L\: k]=d'\endSb N_k(n,L,d'm/d).\tag 17$$
Clearly $N_k(n,L,d'm/d)\le N_L(n,1,d'm/d)$ always. The case where
$m\ge 2g_K-1$ of Corollary 2 follows from Theorem 2, Corollary 1, Lemma 6,
and (15)-(17).
Now suppose $m\le 2g_K-2$. We have the trivial bound
$N_k(n,K,m)\le N_K(n,1,m)$. As remarked following the
statement of Corollary 2, we have
$S_K(n,1)q^{nm}\ll q^{m({n+1\over 2}+\varepsilon )}$ when $m<2g_K-1$ by
Lemmas 2 and 3.
Thus, this case of Corollary 2 follows
directly from Theorem 2.
\enddemo

\head Proof of Theorem 1\endhead

As stated before, we will use Corollary 2 and (2) to prove Theorem 1.
Throughout this section the function field $k$ is fixed and, as before,
we simply write $q$ for $q_k$. We write $e=[k\:\bfq (T)]$ as in
the statement of Theorem 1. All implicit constants depend only on (at most)
$n$, $d$,
$e$, $k$ and $\varepsilon$.

Corollary 2 to Theorem 1
uses the quantity $\delta _n(K/k)$. It turns out simpler
to
use $\delta _2(K/k)$. We thus need to compare the two quantities and
get some useful estimates. But  first we show that these quantities
exist in the first place, since it is not a priori obvious that
there is an $\alpha \in K$ with $k(\alpha )=K$, for example.

\proclaim{Lemma 7} Let $K\supseteq k$ be a function field with
$q_K=q$. Then
there is an $\alpha \in K$ with $k(\alpha )=K$. In particular,
for all $n\ge 2$ there is a point $P\in \bp ^{n-1}(K)$ with $k(P)=K$.
\endproclaim

\demo{Proof} The first assertion follows immediately from the
primitive element theorem (see [H, p. 287], for example) and
the fact, proven in Lemma 6 above, that there are only finitely many
intermediate subfields. As for the second assertion, write $K=k(\alpha )$.
Then $K$ is generated by $(1\:\alpha )\in\bp ^1(K)$ and more generally
$(1\:\alpha \:\cdots \:\alpha )\in\bp ^{n-1}(K)$.
\enddemo

\proclaim{Lemma 8} Let $d>1$ and  let
$K\supseteq k$ be a function field with $d=[K\: k]$ and $q_K=q.$ Then
$$\delta _2(K/k)+1-d\le\delta _i(K/k)\le\delta _j(K/k)$$
for all $2\le j\le i$.
Also
$${g_K\over d-1}-c_1(k,d)\le de\delta _2(K/k)\le g_K+c_2(k,d),$$
where $c_1(k,d)$ and $c_2(k,d)$ are positive integers depending only
on $k$ and $d$. Thus
$${g_K\over d-1}-c_3(k,d)\le de\delta _i(K/k)\le g_K+c_2(k,d)$$
for all $i\ge 2$, where $c_3(k,d)=c_1(k,d)+de(d-1)$.
\endproclaim

\demo{Proof}
Suppose $\delta _j(K/k)=h(P)$ for $P\in\bp ^{j-1}(K)$ with $k(P)=K$.
Without loss of generality we have $P=(1\:\alpha _1\:\cdots \:\alpha _{j-1})$.
But then $P'=(1\: \cdots \: 1\: \alpha _1\:\cdots\:\alpha _{j-1})$
also generates $K$ over $k$ for any number of $1$'s, and clearly
$h(P')=h(P)$. Thus, $\delta _i(K/k)\le\delta _j(K/k)$ whenever $i\ge j\ge 2$.

Now suppose $i>2$ and write
$\delta _i(K/k)=h(1:\alpha _1:\cdots :\alpha _{i-1}),$
where $K=k(\alpha _1,\ldots ,\alpha _{i-1})$. Write $K_s$ for the maximal
separable extension of $k$ contained in $K$ and set $s=[K_s\: k]$.
Arguing exactly as in the proof of [T1, Lemma 7],
we claim that there exist polynomials $f_1,...,f_{i-1}$ in $\bfq [T]$,
either zero or of degree at most $s-1$, such that
$K_s\subseteq k(\alpha)$ for $\alpha=\sum_{l=1}^{i-1} f_l\alpha_l$.
This is trivially true if $s=1$, so assume $s>1$.
For $l=1,\ldots ,s$ let $\sigma_l: K\rightarrow \fptbar$ be the
$k$-homomorphisms of $K$.
By induction on $i$ we easily deduce that for each nonzero homogeneous
polynomial
$P(X_1,...,X_{i-1})\in K[X_1,...,X_{i-1}]$ of degree $A$ there
exist elements $f_1,...,f_{i-1}\in\bfq [T]$, either zero or of degree at
most $A$,
such that $P(f_1,...,f_{i-1})\neq 0$.
We let
$$P(X_1,...,X_n)=\prod_{l=2}^{s}\sum_{j=1}^{i-1}
(\sigma_1(\alpha_j)-\sigma_l(\alpha_j))X_j.$$
Since the $\sigma_l$ are pairwise distinct $k$-homomorphisms on $K$ and
$K=k(\alpha_1,\ldots ,\alpha_{i-1}),$ we conclude
that for each $l>1$ there exists an $\alpha_j$ among
$\alpha_1,...,\alpha_{i-1}$
with $\sigma_1(\alpha_j)\neq \sigma_l(\alpha_j)$.
Therefore $P$ is not the zero polynomial. Furthermore the degree of
$P$ is $s-1$. Hence we can find
$f_1,\ldots ,f_{i-1}\in\bfq [T]$, either zero or of degree at most $s-1$,
with $P(f_1,\ldots ,f_{i-1})\neq 0$. But this means that
$\sigma_1(\alpha)\neq \sigma_l(\alpha)$ for $l=2,\ldots ,s$
where $\alpha=\sum_{j=1}^{i-1}f_j\alpha_j$.
Therefore $K_s\subseteq k(\alpha)$.

With $\alpha$ as above, an easy calculation shows that
$$h(1\:\alpha)\le h(1\:\alpha_1\:\cdots \:\alpha_{i-1})+(s-1)h(1\: T)=
\delta _i(K/k)+s-1.$$
This suffices to prove that $\delta _2(K/k)-d+1\le\delta _i(K/k)$
in the case where $K_s=K$. If $K_s\neq K$, then
$\alpha$ may not generate the entire field $K$, but some $p^r$th root
$\theta$ of $\alpha$ does
by [S, Proposition III.9.2]. In this case we
have $p^rh(1\:\theta )=h(1\:\alpha )\le \delta _i(K/k)+s-1$ and
again $\delta _2(K/k)-d+1\le\delta _i(K/k).$

The upper bound for $\delta _2(K/k)$ is [W2, Theorem 1.1].
The lower bound is [T1, Lemma 6]. (Although separability is a stated
assumption in \S IV of [T1], the proof of Lemma 6 does not use this.)
\enddemo

As noted in the introduction, all of our $\delta _2(K/k)$ (since they are
the height of some point in $\bp ^1(K)$) are necessarily of the form $m/de$ for
some non-negative integer $m$. We will need the following estimate.

\proclaim{Lemma 9} For all $d\ge 2$ we have
$$\sum \Sb [K\: k]=d\\ q_K=q\\ \delta _2(K/k)=m/de\endSb 1\ll q^{(d+1)m}.$$
\endproclaim

\demo{Proof}
Certainly the number of function fields
$K$ with $\delta _2(K/k)=m/de$ for a given $m$ is no greater than the number
of $\alpha$ of degree $d$ over $k$ and $h(1\: \alpha )=m/de$.
Each such $\alpha$ has
a defining polynomial $P$ of degree $d$ over $k$, and by [RT, Lemma 4.9] we
have $h(P)=dh(1\:\alpha )=m/e$, where $h(P)$ denotes the height of the
coefficient
vector of $P$, which we view as a point in $\bp ^d(k)$.
We conclude that the number of $\alpha$ of degree $d$
over $k$ with height $h(1\:\alpha )=m/de$ is no more than $dN_k(d+1,1,m)$.
Thus, by Corollary 1 we have
$$\sum \Sb [K\: k]=d\\ q_K=q\\ \delta _2(K/k)=m/de\endSb 1\le
dN_k(d+1,1,m)\ll q^{(d+1)m}.$$
\enddemo

\proclaim{Lemma 10} For all integers $d>1$ and $m\ge 0$,
$$\sum \Sb [K\: k]=d\\ q_K=q\\ g_K= m\endSb 1\le
\sum \Sb [K\: k]=d\\ q_K=q\\ g_K\le m\endSb 1\ll q^{(d+1)m}.$$
\endproclaim

\demo{Proof} By Lemmas 8 and 9
$$\aligned
\sum \Sb [K\: k]=d\\ q_K=q\\ g_K\le m\endSb 1&\le
\sum \Sb [K\: k]=d\\ q_K=q\\ de\delta _2(K/k)\le m+c_2(k,d)\endSb 1\\
&=\sum _{i=0}^{m+c_2(k,d)}
\sum \Sb [K\: k]=d\\ q_K=q\\ de\delta _2(K/k)=i\endSb 1\\
&\ll \sum _{i=0}^{m+c_2(k,d)}q^{(d+1)i}\\
&\ll q^{(d+1)m}.\endaligned$$
\enddemo

\demo{Proof of Theorem 1}
We assume that $\varepsilon >0$
and that $n$ and $d$ are positive integers with $n\ge 4$ and $d>1$ initially.
By (2) and Corollary 2 to Theorem 2,
$$\aligned N_k(n,d,m)&=
\sum \Sb [K\: k]=d\\ q_K=q\endSb N_k(n,K,m)\\
&=
\sum \Sb [K\: k]=d\\ q_K=q\\ de\delta _n(K/k)\le m\endSb
S_K(n,1)q_K^{nm}
+O\left (\sum \Sb [K\: k]=d\\ q_K=q\\ 2g_K-1\le m\endSb
q^{nm/2}\right )\\
&\qquad
+O\left (\sum \Sb [K\: k]=d\\ q_K=q\\ de\delta _n(K/k)\le m< 2g_K-1\endSb
q^{m({n+1+\varepsilon\over 2})}\right ).
\endaligned\tag 18$$

First, for the main term we claim that
$$S_k(n,d)=\sum \Sb [K\: k]=d\\ q_K=q\endSb S_K(n,1)$$
converges whenever $n>d+2$. Moreover, we claim that
$$\sum \Sb [K\: k]=d\\ q_K=q\\ de\delta _2(K/k)\ge m\endSb S_K(n,1)
\ll {q^{m(d+2+\varepsilon -n)}\over 1-q^{d+2+\varepsilon -n}}\tag 19$$
whenever $n>d+2+\varepsilon$.
Indeed, by Lemmas 2 and 3 $S_K(n,1)\ll q^{g_K(1+\varepsilon -n)}$ and
therefore by Lemmas 8 and 9
$$\aligned \sum \Sb [K\: k]=d\\ q_K=q\\ de\delta _n(K/k)\ge m\endSb S_K(n,1)
&\ll
\sum \Sb [K\: k]=d\\ q_K=q\\ de\delta _n(K/k)\ge m\endSb
q^{g_K(1+\varepsilon -n)}\\
&\le
\sum \Sb [K\: k]=d\\ q_K=q\\ de\delta _2(K/k)\ge m\endSb
q^{g_K(1+\varepsilon -n)}\\
&\ll
\sum \Sb [K\: k]=d\\ q_K=q\\ de\delta _2(K/k)\ge m\endSb
q^{de\delta _2(K/k)(1+\varepsilon -n)}\\
&=\sum _{i=m}^{\infty}
\sum \Sb [K\: k]=d\\ q_K=q\\ de\delta _2(K/k)=i\endSb
q^{i(1+\varepsilon -n)}\\
&\ll
\sum _{i=m}^{\infty}q^{i(d+2+\varepsilon -n)}\\
&\ll
{q^{m(d+2+\varepsilon -n)}\over 1-q^{d+2+\varepsilon -n}}
,\endaligned$$
proving (19) and also that the sum defining $S_k(n,d)$ converges.

Now for the error terms. By Lemma 10
$$\sum \Sb [K\: k]=d\\ q_K=q\\ 2g_K-1\le m\endSb q^{nm/2}
\ll q^{nm/2}q^{m(d+1)/2}=q^{{m\over 2}(n+d+1)},\tag 20$$
and by Lemmas 8 and 9
$$\aligned \sum \Sb [K\: k]=d\\ q_K=q\\ de\delta _n(K/k)\le m\endSb
q^{m({n+1+\varepsilon\over 2})}&\le
\sum \Sb [K\: k]=d\\ q_K=q\\ de\delta _2(K/k)\le m+d-1\endSb
q^{m({n+1+\varepsilon\over 2})}\\
&=\sum _{i=0}^{m+d-1}
\sum \Sb [K\: k]=d\\ q_K=q\\ de\delta _2(K/k)=i\endSb
q^{m({n+1+\varepsilon\over 2})}\\
&\ll \sum _{i=0}^{m+d-1}q^{m({n+1+\varepsilon\over 2})}q^{i(d+1)}\\
&\ll q^{m({n+1+\varepsilon\over 2})}q^{m(d+1)}\\
&=q^{{m\over 2}(n+2d+3+\varepsilon )}.\endaligned\tag 21$$
The proof of Theorem 1 is completed by (18)-(21).
\enddemo

We now turn to the Corollary. At this point we are forced to distinguish
between separable and inseparable extensions of the field $k$. Similar
to (2), we write
$$\nksep (n,d,m)=\sum ^{\text{sep}}\Sb [K\: k]=d\\ q_K=q_k\endSb N_k(n,K,m)$$
where the superscript on the summation indicates that we sum only over
separable extensions $K$. Similarly, we write $\nfksep (n,d,m)$ for
the number of forms counted in $NF_k(n,d,m)$ where each $k\big (P(\uL _i)
\big )$ is a separable extension of $k$.

We first show that
$$\nfksep (n,d,m)={1\over d}\nksep (n,d,m).\tag 22$$
Towards that end, suppose $F(\uX )=\prod _{i=1}^dL_i(\uX )$
is a decomposable form counted
in $\nfksep (n,d,m)$. Then each $P(\uL _i)$ generates a separable extension
of $k$ of degree $d$ and effective degree $d$. Further, by the separability
assumption these points are pairwise distinct.
From [RT, Lemma 4.9] we have $h\big (P(F)\big )=dh\big (P(\uL _i)\big )$ for
each $i=1,\ldots ,d$, so that each $P(\uL _i)$ is counted in $\nksep (n,d,m)$
(recall that the counting function $N_k$ takes into account the effective
degree).
Conversely, suppose $P\in\bp ^{n-1}(\fptbar )$ is counted in $\nksep (n,d,m)$.
Then
$P=P(\uL _1)$ for an
$\uL _1\in k(P)^n$ with $d$ pairwise distinct conjugates $\uL _1,\ldots ,
\uL _d$ and
$F(\uX )=\prod _{i=1}^dL _i(\uX )\in k[\uX ]$ is counted in $\nfksep (n,m,d)$.
Since each $P(\uL _i)$ here is counted in $\nksep (n,m,d)$, this proves (22).

As in the proof of Lemmas 6 and 8, for an extension $K$ of $k$ we write
$K_s$ for the maximal separable extension of $k$ contained in $K$. For
such a field $K$ we have $[K\: K_s]=p^r$ for some integer $r\ge 0$. As
remarked above in the proofs of Lemmas 6 and 8,
if $P=(\alpha _0\:\cdots \: \alpha _n)
\in\bp ^n(K)$ with $k(P)=K$, then $Q=
(\alpha _0^{p^r}\:\cdots \: \alpha _n
^{p^r})\in\bp ^n(K_s)$ with $K_s=k(Q)$ and $p^rh(P)=h(Q)$.
Hence $N_k(n,K_s,m)=N_k(n,K,m)$ (recall that the definition
of $N_k$ takes into account the effective degree of the extension), so that
$$\gathered
N_k(n,d,m)=\sum \Sb p^r|d\endSb \nksep (n,d/p^r,m)\\
NF_k(n,d,m)=\sum \Sb p^r|d\endSb \nfksep (n,d/p^r,m).\endgathered\tag 23$$

We claim that
$$\nksep (n,d,m)=\cases N_k(n,d,m)-N_k(n,d/p,m)&\text{if $p|d$,}\\
N_k(n,d,m)&\text{if $p\nmid d$.}\endcases\tag 24$$
We prove this by induction on the highest power of $p$ dividing $d$. This
is clearly true if $p\nmid d$, so assume $p^r$ is the highest power of
$p$ dividing $d$ with $r>0$. Then by (23) and the induction hypothesis,
$$\aligned
N_k(n,d,m)&=\sum _{i=0}^r\nksep (n,d/p^i,m)\\
&=\nksep (n,d,m)+\sum _{i=1}^r\nksep (n,d/p^i,m)\\
&=\nksep (n,d,m)+N_k(n,d/p^r,m)+
\sum _{i=1}^{r-1}N_k(n,d/p^i,m)-N_k(n,d/p^{i+1},m)\\
&=\nksep (n,d,m)+N_k(n,d/p,m).\endaligned$$
The proof of the Corollary is completed with $(22)-(24)$ and Theorem 1.

Finally, we turn to our remark regarding possible asymptotic results.
Though we did not need it for the proof of Lemma
9 above, it is known (see [T3, Theorem 1])
that $N_k(2,d,m)$ is actually asymptotic
to $dN_k(d+1,1,m)$. In particular,
$$N_k(2,d,m)\gg\ll q^{m(d+1)}.$$
On the other hand, we clearly have $N_k(n,d,m)\ge N_k(2,d,m)$ for all
$n\ge 2$. This immediately implies that $N_k(n,d,m)$ cannot be asymptotic
to $cq^{nm}$ for any real $c$ when $n<d+1$.

\head Acknowledgements\endhead
The second author thanks Paula Tretkoff for inspiring questions and the
Institut des Hautes \'Etudes Scientifiques in Bures-sur-Yvette, where parts
of this work have been done, for the wonderful hospitality.

\enddocument